\documentclass[11pt]{article}
\usepackage{amsmath,amssymb,theorem, fullpage}
\usepackage{graphicx}
\usepackage{subfigure}
\usepackage{psfrag}
\usepackage{color}
\usepackage{enumerate}
\usepackage{epstopdf}
\usepackage[ruled, lined, linesnumbered]{algorithm2e}

\newtheorem{thm}{Theorem}[section]
\newtheorem{conj}[thm]{Conjecture}
\newtheorem{lem}[thm]{Lemma}
\newtheorem{ques}[thm]{Question}
\theorembodyfont{\rmfamily}
\def\pf{\bigskip\noindent {\bf Proof.}~~}

\def\dfn#1{{\sl #1}}

\def\pf{\bigskip\noindent {\bf{Proof.}}~~}

\usepackage[natural]{xcolor}
\usepackage{listings}
\title{Strong list-chromatic index of subcubic graphs}

\author{Tianjiao Dai$^1$\thanks{Email: dtjmath@163.com.}, \ \ Guanghui Wang$^1$\thanks{Corresponding author. Email: ghwang@sdu.edu.cn. Supported by the National Natural Science Foundation of China (11471193,11631014), the Foundation for Distinguished Young Scholars of Shandong Province (JQ201501).}, \ \ Donglei Yang$^1$\thanks{Email: yangdonglei\_sdu@163.com.}, \ \ Gexin Yu$^{2,3}$\thanks{Email: gyu@wm.edu. Supported by the NSA (H98230-16-1-0316) and  NSFC (11728102).}\\
$^{1}$\small Department of Mathematics, Shandong University, Jinan, Shandong, China.\\
$^2$\small Department of Mathematics, The College of William and Mary, Williamsburg, VA, USA.\\
$^{3}$\small Department of Mathematics, Central China Normal University, Wuhan, Hubei, China.\\
}
\date{}

\begin{document}
\maketitle


\begin{abstract}
A strong $k$-edge-coloring of a graph G is an edge-coloring with $k$ colors in which every color class is an induced matching. The strong chromatic index of $G$, denoted by $\chi'_{s}(G)$, is the minimum $k$ for which $G$ has a strong $k$-edge-coloring.
 In 1985, Erd\H{o}s and Ne\v{s}et\v{r}il conjectured that $\chi'_{s}(G)\leq\frac{5}{4}\Delta(G)^2$, where $\Delta(G)$ is the maximum degree of $G$. When $G$ is a graph with maximum degree at most 3, the conjecture was verified independently by Andersen and Hor\'{a}k, Qing, and Trotter. In this paper, we consider the list version of strong edge-coloring. In particular, we show that every subcubic graph has strong list-chromatic index at most 11 and every planar subcubic graph has strong list-chromatic index at most 10.

\end{abstract}
\bigskip
\noindent {\textbf{Keywords}: Subcubic graphs; Strong list-chromatic index; Combinatorial Nullstellensatz

\baselineskip 18pt
\section{Introduction}
All graphs in this paper are finite and simple. A \dfn{strong $k$-edge-coloring} of a graph $G$ is a coloring $\phi: E(G)\longrightarrow [k]$ such that if any two edges $e_1$ and $e_2$ are either adjacent to each other or adjacent to a common edge, then $\phi(e_1)\neq \phi(e_2)$. In other words, the edges in each color class give an induced matching in the graph;  that is, any two vertices belonging to distinct edges with the same color are not adjacent. The \dfn{strong chromatic index} of $G$, denoted by $\chi'_{s}(G)$, is the minimum $k$ for which $G$ has a strong $k$-edge-coloring.

The following conjecture was proposed by Erd\H{o}s and Ne\v{s}et\v{r}il in 1985 at Prague.

\begin{conj}\label{conj1}\emph{~\cite{Erdos1,Erdos2}}
If $G$ is a graph with maximum degree $\Delta(G)$, then
\begin{equation*}
\chi'_{s}(G)\leq
\begin{cases}
  \frac{5}{4}\Delta(G)^2, & \text{if $\Delta(G)$ is even,}\\
  \frac{5}{4}\Delta(G)^2-\frac{1}{2}\Delta(G)+\frac{1}{4}, & \text{if $\Delta(G)$ is odd.}\end{cases}
  \end{equation*}
\end{conj}

Note that there are examples showing that the conjectured upper bound is tight (i.e. blow-ups of a $5$-cycle). Andersen \cite{And} and independently Hor\'{a}k, Qing, and Trotter \cite{HQ} showed that $\chi'_{s}(G)\leq 10$ for any graph $G$ with $\Delta(G)=3$, thus settling the first nontrivial case of Conjecture \ref{conj1}.  Cranston~\cite{C} gave an algorithm that uses at most 22 colors for every graph with $\Delta(G)=4$, which was improved to $21$ very recently by Huang, Santana and Yu~\cite{Y}.    When $\Delta(G)$ is sufficiently large, Molloy and Reed~\cite{MR} proved that $\chi'_{s}(G)\leq 1.998\Delta(G)^2$. Henning Bruhn and Felix Joos~\cite{HB} showed that $\chi'_{s}(G)\leq 1.93\Delta(G)^2$. Recently, Bonamy, Perrett, and Postle~\cite{BPP} improved the upper bound to $1.835\Delta(G)^2$.

In this article, we study the list version of strong edge-coloring. For each $e\in E(G)$, let $L(e)$ be the list of available colors of $e$, and let $L=\{L(e): e\in E(G)\}$.  The graph $G$ is {\em strongly $L$-edge-colorable} if there exists a strong edge coloring $c$ of $G$ such that $c(e)\in L(e)$ for every $e\in E(G)$. For a positive integer $k$, a graph $G$ is {\em strongly $k$-edge-choosable} if $G$ is strongly $L$-edge colorable for every $L$ with $|L(e)|\ge k$ for all $e\in E(G)$. The {\em strong list-chromatic index}, denoted by $\chi'_{s,l}(G)$, is the minimum positive integer $k$ for which $G$ is strongly $k$-edge-choosable. Note that $\chi'_{s}(G)\leq\chi'_{s,l}(G)$ for every graph $G$.

The probablistic arguments that Molloy-Reed and Bonamy-Perrett-Postle used to give upper bounds of $\chi'_s$ on graphs of large $\Delta(G)$ actually also work for the strong list-chromatic index. So we have $\chi'_{s,l}(G)\le 1.835\Delta(G)^2$ for large $\Delta(G)$.  Ma, Miao, Zhu, Zhang and Luo~\cite{ML} proved that the strong list-chromatic index of a subcubic graph with maximum average degree less than $\frac{15}{7}, \frac{27}{11}, \frac{13}{5}, \frac{36}{13}$ is at most $6, 7, 8, 9$, respectively.  More results of this kind can be found in \cite{Z}.

In this paper, we prove the following result.

\begin{thm}\label{th1}
If $\Delta(G)\le 3$, then $\chi'_{s,l}(G)\leq 11$.
\end{thm}

For planar graphs, we actually can do a little better.

\begin{thm}\label{planar}
If $G$ is a subcubic planar graph, then $\chi'_{s,l}(G)\leq 10$.
\end{thm}

Recall that Andersen~\cite{And} and Hor\'{a}k, Qing, and Trotter~\cite{HQ} proved that $\chi'_s(G)\le 10$ if $\Delta(G)\le 3$. Kostochka et. al.~\cite{K} proved that $\chi'_s(G)\le 9$ under the additional assumption that $G$ is planar. We do not feel that our results are optimal, but it may involve substantial work to improve them.

One of the main tools we use is Hall's Theorem.

\begin{lem}\label{lemma1} \emph{(Hall \cite{Ha})} Let $A_{1},...,A_{n}$ be $n$ subsets of a set U. Distinct representatives of $\{A_{1},...,A_{n}\}$ exist if and only if for all $k$, $1\leq k\leq n$ and every choice of subcollection of size $k$, $\{A_{i_{1}},...,A_{i_{k}}\}$, we have $|A_{i_{1}}\bigcup ...\bigcup A_{i_{k}}|\geq k$.
\end{lem}

Another tool we use is the Combinatorial Nullstellensatz.

\begin{lem}\label{lemma2} \emph{(Alon \cite{Alon}, Combinatorial Nullstellensatz)} Let $\mathbb{F}$ be an arbitrary field, and let
$P = P(x_1, x_2, \ldots, x_n)$ be a polynomial in $\mathbb{F}[x_1, x_2,\ldots, x_n]$. Suppose that the degree $deg(P)$ of $P$ equals $\sum^n_{i=1}k_i$, where each $k_i$ is a non-negative integer, and the coefficient of $x_1^{k_1}x_2^{k_2}\cdots x_n^{k_n}$ in $P$ is non-zero. Then if $S_1, S_2,\ldots, S_n$ are subsets of $\mathbb{F}$ with $|S_i| > k_i, i=1,2,\ldots,n$, there exist $s_1\in S_1,s_2\in S_2,\ldots, s_n\in S_n$
so that $P(s_1, s_2, \ldots, s_n)\neq0$.
\end{lem}

We use MATLAB to calculate the coefficients of specific monomials. Let $P=P(x_1, x_2, \ldots, x_n)$ be a polynomial in $n$ variables, where $n\geq 1$. By $c_p(x^{k_1}_{1}x^{k_2}_{2}\cdots x^{k_n}_{n})$, we denote the coefficient of the monomial $x^{k_1}_{1}x^{k_2}_{2}\cdots x^{k_n}_{n}$ in $P$, where $k_i$ $(1\leq i\leq n)$ is a non-negative integer. The codes are listed in the Appendix.

\section{Basic properties}
 Consider $(G,L)$ such that $G$ is not $L$-choosable but any proper subgraph of $G$ is $L$-choosable.  Clearly,  $G$ is connected. In this section, we will show that if $|L(e)|\ge 10$ for each $e\in E(G)$, then $G$ is cubic and has no cycles of length at most five.


We first introduce some notation. An \dfn{$i$-vertex} is a vertex of degree $i$ in our graphs. An \dfn{$i$-cycles} is a cycle of length $i$ in graphs. A \dfn{partial coloring} of $G$ is a coloring of a proper subgraph of $G$. Given edges $e$ and $ e'$ in $G$, we say that {\em $e$ sees $e'$} if either $e$ and $e'$ are adjacent, or there is another edge $e''$ adjacent to both $e$ and $e'$. Note that even if $e$ sees $e'$ in G, $e$ does not necessarily see $e'$ in a proper subgraph of $G$. Additionally, we will also say that {\em $e$ sees a color $\alpha$} if $e$ sees an edge $e'$ of color $\alpha$. Let $\phi$ be a partial coloring of $G$. For $e\in E(G)$, let $C_{\phi}(e)$ denote the set of colors seen by $e$, and let $A_\phi(e)=L(e)\setminus C_{\phi}(e)$. For $v\in V(G)$, $H\subseteq G$, let $d(v,H)$ with respect to $v$ be the minimum of the lengths of the $u$-$v$ paths of $G$ where $u\in V(H)$. 

\begin{lem}\label{lemma3.1}
$G$ is cubic.
\end{lem}

\pf By way of contradiction, we assume that $d(v)\le 2$ for some $v\in V(G)$. By the minimality of $G$, $G-v$ has an $L$-coloring $\phi$.    First let $v$ be a $1$-vertex incident with the edge $e$. Since $|C_{\phi}(e)|\le 6$, $|A_{\phi}(e)|\ge 4$, so $e$ can be colored. Let $v$ be a $2$-vertex with incident edges $e_1$ and $e_2$.  Since $|C_{\phi}(e_i)|\le 8$ for $i=1,2$, $|A_{\phi}(e_i)|\ge 2$. So we can color $e_1$ and $e_2$ in any order. $\Box$

\begin{lem}\label{lemma3.2}
$G$ has no triangles.
\end{lem}

\pf Suppose that $G$ contains a triangle: $v_{1}v_{2}v_{3}v_{1}$ (see Fig.1 (1)). By the minimality of $G$, let $\phi$ be an $L$-coloring of the subgraph $H=G-v_{1}$. Note that $|A_\phi(e_i)|\geq3$, for $i=1,2$ and $|A_\phi(e_3)|\geq1$. Then $\phi$ can be extended to an $L$-coloring of $G$ by Lemma~\ref{lemma1}. $\Box$

\begin{figure}[htbp]
\centering
\includegraphics[scale=0.45]{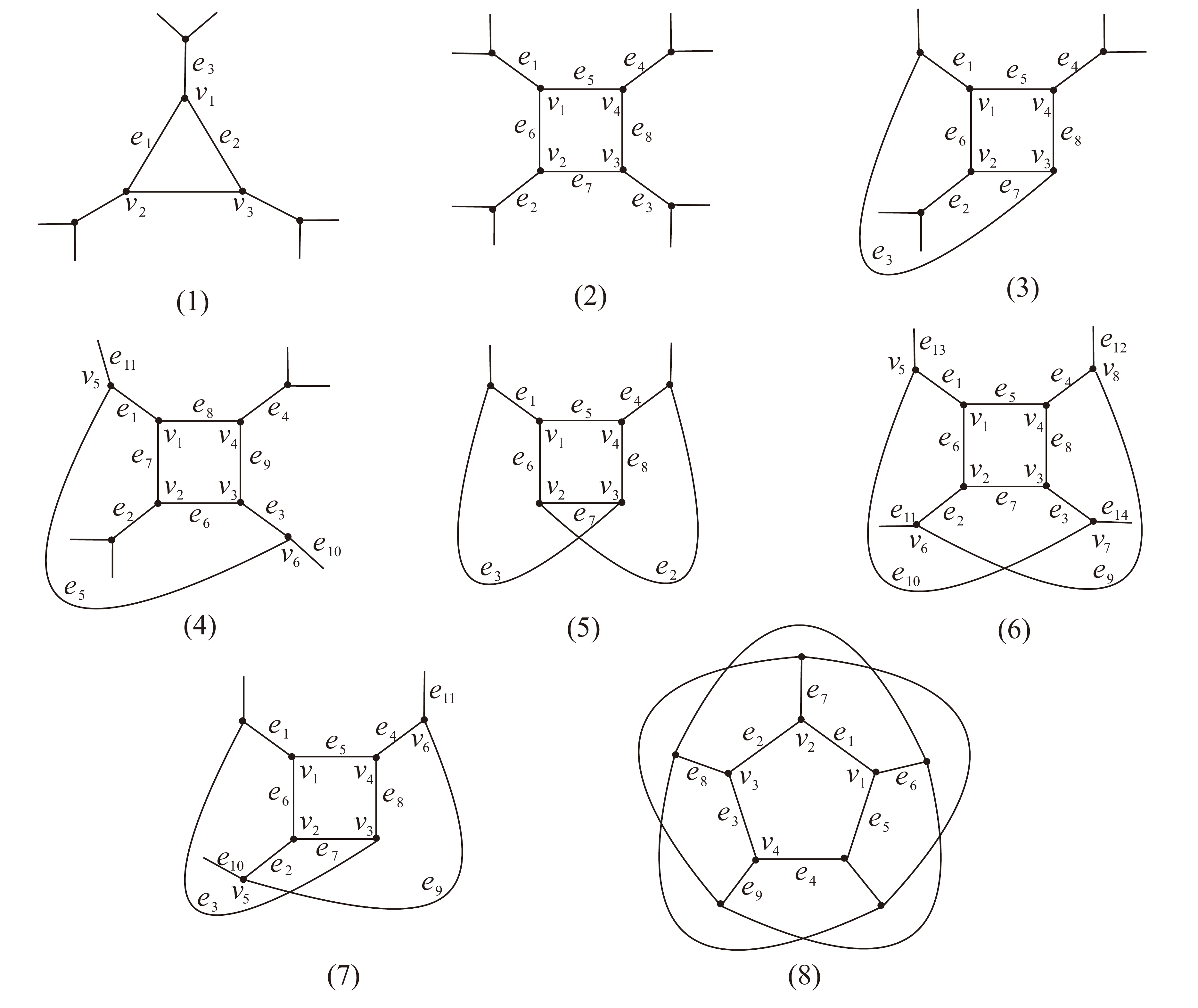}
\caption{Structures with $3$-, $4$-, and $5$-cycles}
\label{C1}
\end{figure}

\begin{lem}\label{lemma3.3}
$G$ has no 4-cycles.
\end{lem}

\pf Suppose that $G$ contains a $4$-cycle. By Lemma~\ref{lemma3.2}, each $4$-cycle must be an induced $4$-cycle, and we divide all $4$-cycles into three classes. The possible local structures about a 4-cycle are shown in Figure 1 (2)-(7). For $2\le i\le 7$, let $H_i$ be the subgraph of $G$ obtained by removing the vertices with labels in Figure~\ref{C1} (i) for $2\le i\le 7$. By minimality of $G$, $H_i$ has an $L$-coloring $\phi$. Let $A_i=A_\phi(e_i)$ for each labelled edge $e_i$ in the figures.

When $e_1$ does not see $e_3$ and $e_2$ does not see $e_4$, we consider $H_2$.

In $H_2$, we have $|A_i|\geq4$ for $i=1,2,3,4$ and $|A_i|\geq6$ for $i=5,6,7,8$. By Lemma~\ref{lemma1}, we may assume that for some $I\subseteq [8]$, $|\bigcup_{i\in I} A_i|<|I|$. So $|I|>6$, and $|I|\in \{7,8\}$. By symmetry, let $1,3\in I$. Then $|A_1\cup A_3|<8$, so there exists $\alpha\in A_1\cap A_3$, and we color $e_1$ and $e_3$ with $\alpha$. Let $A_i'=A_i-\{\alpha\}$. Now for  $J\subseteq [8]-\{1,3\}$, if $|\bigcup_{i\in J} A_i'|<|J|$, then $|J|>3$, which implies that $J\cap \{5,6,7,8\}\not=\emptyset$, so $|J|\ge 6$.  Then $|J|=6$,  $|A_i'|\geq5$ for $i\in \{5,6,7,8\}$ and $A_2'\cap A_4'\not=\emptyset$.  Color $e_2$ and $e_4$ with $\beta\in A_2'\cap A_4'$, then we can color $e_5, e_6, e_7, e_8$ in any order, a contradiction.

When $e_1$ sees $e_3$ and $e_2$ does not see $e_4$, we consider $H_3$ and $H_4$.

Consider $H_3$. Note that $|A_i|\geq7$ for $i=1,3,5,6,7,8$ and $|A_i|\geq4$ for $i=2,4$.  By Lemma~\ref{lemma1}, we may assume that for some $I\subseteq [8]$, $|\bigcup_{i\in I} A_i|<|I|$. Then $I=[8]$,  $|A_i|\geq7$ for $i\in \{1,3,5,6,7,8\}$ and $A_2\cap A_4\not=\emptyset$. Color $e_2$ and $e_4$ with $\alpha\in A_2\cap A_4$, and then we can color the rest of edges one by one, a contradiction.

Consider $H_4$. Note that $|A_i|\geq4$ for $i=2,4,10,11$, $|A_5|\geq6$ and $|A_i|\geq8$ for $i\in \{1,3,6,7,8,9\}$. By Lemma~\ref{lemma1}, we may assume that for some $I\subseteq [11]$, $|\bigcup_{i\in I} A_i|<|I|$.  Then $|I\cap \{1,3,6,7,8,9\}\not=\emptyset$, so $|I|\ge 9$.

We consider the following cases.

Case 1:  $\{6,11\}\subset I$ (or by symmetry, $\{7,10\}\subset I$).  Then $A_6\cap A_{11}\not=\emptyset$ and color $e_6$ and $e_{11}$ with $\alpha\in A_6\cap A_{11}$.  For $i\in [11]-\{6,11\}$, let $A_i'=A_i-\{\alpha\}$. Then for some $J\subseteq [11]-\{6,11\}$, $|\bigcup_{i\in J} A_i'|<|J|$. Then $|J|\ge 8$. So at most one of $2,4,5,7,10$ is not in $J$.

\begin{itemize}
\item $5\in J$. We may assume that $2\in J$ as well (or by symmetry, $4\in J$).  Then $A_2'\cap A_5'\not=\emptyset$, and we color $e_2$ and $e_5$ with $\beta\in A_2'\cap A_5'$. Let $A_i''=A_i-\{\alpha, \beta\}$ for $i\in [11]-\{2,5,6,11\}$. Then for some $K\subseteq [11]-\{2,5,6,11\}$, $|\bigcup_{i\in K} A_i''|<|K|$.  So $K\cap \{1,3,7,8,9\}\not=\emptyset$ and thus $|K|=7$.  Now we can color $e_7$ and $e_{10}$ with $\gamma\in A_7''\cap A_{10}''$, and by Lemma~\ref{lemma1}, color the rest of the edges.

\item $5\not\in J$.  Then $|J|=8$, so $A_2'\cap A_4'\not=\emptyset$.  Color $e_2$ and $e_4$ with $\beta\in A_2'\cap A_4'$. Let $A_i''=A_i-\{\alpha, \beta\}$ for $i\in [11]-\{6,11, 2,4\}$.  Then  for some $K\subseteq [11]-\{2,4,6,11\}$, $|\bigcup_{i\in K} A_i''|<|K|$.  So $K\cap \{1,3,7,8,9\}\not=\emptyset$ and thus $|K|=7$.  Now we can color $e_7$ and $e_{10}$ with $\gamma\in A_7''\cap A_{10}''$, and by Lemma~\ref{lemma1}, color the rest of the edges.
\end{itemize}

Case 2: $\{6,11\}\not\subset I$ and $\{7,10\}\not\subset I$. Then $|I|=9$ and $2,4,5\in I$. So $A_2\cap A_5\not=\emptyset$. Color $e_2$ and $e_5$ with $\alpha\in A_2\cap A_5$.   For $i\in [11]-\{2,5\}$, let $A_i'=A_i-\{\alpha\}$. Then for some $J\subseteq [11]-\{2,5\}$, $|\bigcup_{i\in J} A_i'|<|J|$. Then $|J|\ge 8$. Then $\{6,11\}\subset J$ (or by symmetry $\{7,10\}\subset J$), and $A_6\cap A_{11}\not=\emptyset$ and color $e_6$ and $e_{11}$ with $\beta\in A_6\cap A_{11}$.  Let $A_i''=A_i-\{\alpha, \beta\}$ for $i\in [11]-\{2,5,6,11\}$. Then  for some $K\subseteq [11]-\{2,5,6,11\}$, $|\bigcup_{i\in K} A_i''|<|K|$. So $K\cap \{1,3,7,8,9\}\not=\emptyset$ and thus $|K|=7$. Now we can color $e_7$ and $e_{10}$ with $\gamma\in A_7''\cap A_{10}''$, and by Lemma~\ref{lemma1}, color the rest of the edges.

When $e_1$ sees $e_3$ and $e_2$ sees $e_4$, we consider $H_5$, $H_6$ and $H_7$.

Consider $H_5$. Note that $|A_i|\geq7$ for $i=1,2,3,4$ and $|A_i|\geq8$ for $i=5,6,7,8$. By Lemma~\ref{lemma1}, we may assume that for some $I\subseteq [8]$, $|\bigcup_{i\in I} A_i|<|I|$. Clearly, no such $I$ exists, a contradiction.

Consider $H_6$. First, $|A_\phi(e_i)|\geq4$ for $i=11,12,13,14$. We can make it that the colors of $e_i$ are different by Lemma 1.4 for $i=11,12,13,14$. Then we note that $|A_i|\geq6$ for $i=1,2,3,4$, $|A_i|\geq8$ for $i=5,6,7,8$ and $|A_i|\geq4$ for $i=9,10$. By Lemma~\ref{lemma1}, we may assume that for some $I\subseteq [10]$, $|\bigcup_{i\in I} A_i|<|I|$. Then $|I|>8$, so $|I|\in \{9,10\}$. By symmetry, we may assume that $\{4,10\}\subset I$. Then $A_4\cap A_{10}\not=\emptyset$, so color $e_4, e_{10}$ with $\alpha\in A_4\cap A_{10}$. Let $A_i'=A_i-\{\alpha\}$ for $i\in [10]-\{4,10\}$.  Then for some $J\subseteq [10]-\{4,10\}$, $|\bigcup_{i\in J} A_i'|<|J|$.  It implies that $J\cap \{5,6,7,8\}\not=\emptyset$, thus $J=[10]-\{4,10\}$, and $|A_i'|\geq7$ for $i\in \{5,6,7,8\}$ and $A_1'\cap A_9'\not=\emptyset$.  Color $e_1$ and $e_9$ with $\beta\in A_1'\cap A_9'$, then we can color the rest of the edges one by one, a contradiction.

Consider $H_7$. First, $|A_\phi(e_i)|\geq 4$ for $i=10,11$. We can make it that the color of $e_{10}$ is different from the color of $e_{11}$. Then we note that $|A_i|\geq7$ for $i=1,3$, $|A_i|\geq6$ for $i=2,4$, $|A_i|\geq8$ for $i=5,6,7,8$ and $|A_\phi(e_9)|\geq4$. By Lemma~\ref{lemma1}, we may assume that for some $I\subseteq [9]$, $|\bigcup_{i\in I} A_i|<|I|$. Then $I=[9]$ and $A_3\cap A_9\not=\emptyset$. Color $e_3$ and $e_9$ with $\alpha\in A_3\cap A_9$, and then we can color the rest of edges one by one, a contradiction. $\Box$

\begin{lem}\label{lemma3.4}
$G$ has no 5-cycles.
\end{lem}

\pf Suppose that $G$ contains the 5-cycle (see Figure~\ref{C1} (8)). Then by the minimality of $G$, there is an $L$-coloring $\phi$ of  $H=G-\{v_i: i\in [4]\}$.   We want to color $e_i$ with a color $s_i\in A_{\phi}(e_i)$ for $i\in [9]$ such that close ones do not see each other. So we need to find  $s_i\in A_{\phi}(e_{i})$ for $i\in [9]$ such that $P(s_1,s_2,s_3,s_4,s_5,s_6,s_7,s_8,s_9)\neq0$, where
\begin{align*}
P(x_1,x_2,x_3,x_4,x_5,x_6,x_7,x_8,x_9)=\frac{\prod\limits_{1\leq k<l\leq9}(x_k-x_l)}{(x_1-x_9)(x_5-x_8)(x_3-x_6)(x_4-x_7)}.
\end{align*}

Note that $\deg(P)=32$, $|A_\phi(e_2)|\geq6$ and $|A_{\phi}(e_i)|\geq5$ for $i\in [9]-\{2\}$.  Our MATLAB codes show that $c_P(x_1^{4}x_2^{5}x_3^{4}x_4^{4}x_5^{3}x_6^{3}x_7^{3}x_8^{3}x_9^{3})=-6.$ 
By Lemma~\ref{lemma2}, there exist $s_i\in A_i$ for $i\in [9]$ such that $P(s_1,s_2,s_3,s_4,s_5,s_6,s_7,s_8,s_9)\neq0$. Note that the polynomial $P'$ of any other $5$-cycle in $G$ is a subpolynomial of $P$, then $P\not=0$ implies that $P'\not=0$ as well. $\Box$

\medskip

Now we are ready to prove Theorem~\ref{planar}.

\textbf{Proof of Theorem~\ref{planar}}. Let $G$ be a minimal counterexample. By Lemma 2.1-2.4,  the girth of $G$ is at least six. By Euler's formula, $\sum_{v\in V(G)}(2d(v)-6)+\sum_{f\in F(G)}(d(f)-6)=-12.$ It follows that the minimum degree of $G$ is at most two, a contradiction to Lemma~\ref{lemma3.1} that $G$ is $3$-regular. $\Box$

\section{Proof of Theorem~\ref{th1}}

In this section, we give a proof of Theorem~\ref{th1}. Let $(G, L)$ be a minimal counterexample, where $|L(e)|\ge 11$ for each $e\in E(G)$. Without loss of generality, we assume $|L(e)|=11$. By Lemma~\ref{lemma3.1}-~\ref{lemma3.4}, $G$ is $3$-regular and the girth of $G$ is at least six.  Let $v\in V(G)$ with $N(v)=\{v_1,v_2,v_3\}$, and let $N(v_i)-\{v\}=\{w_i, w_i'\}$ for $i\in [3]$. Then $w_1, w_2, w_3$ form an independent set.

\begin{figure}[htbp]
\centering
\includegraphics[scale=0.5]{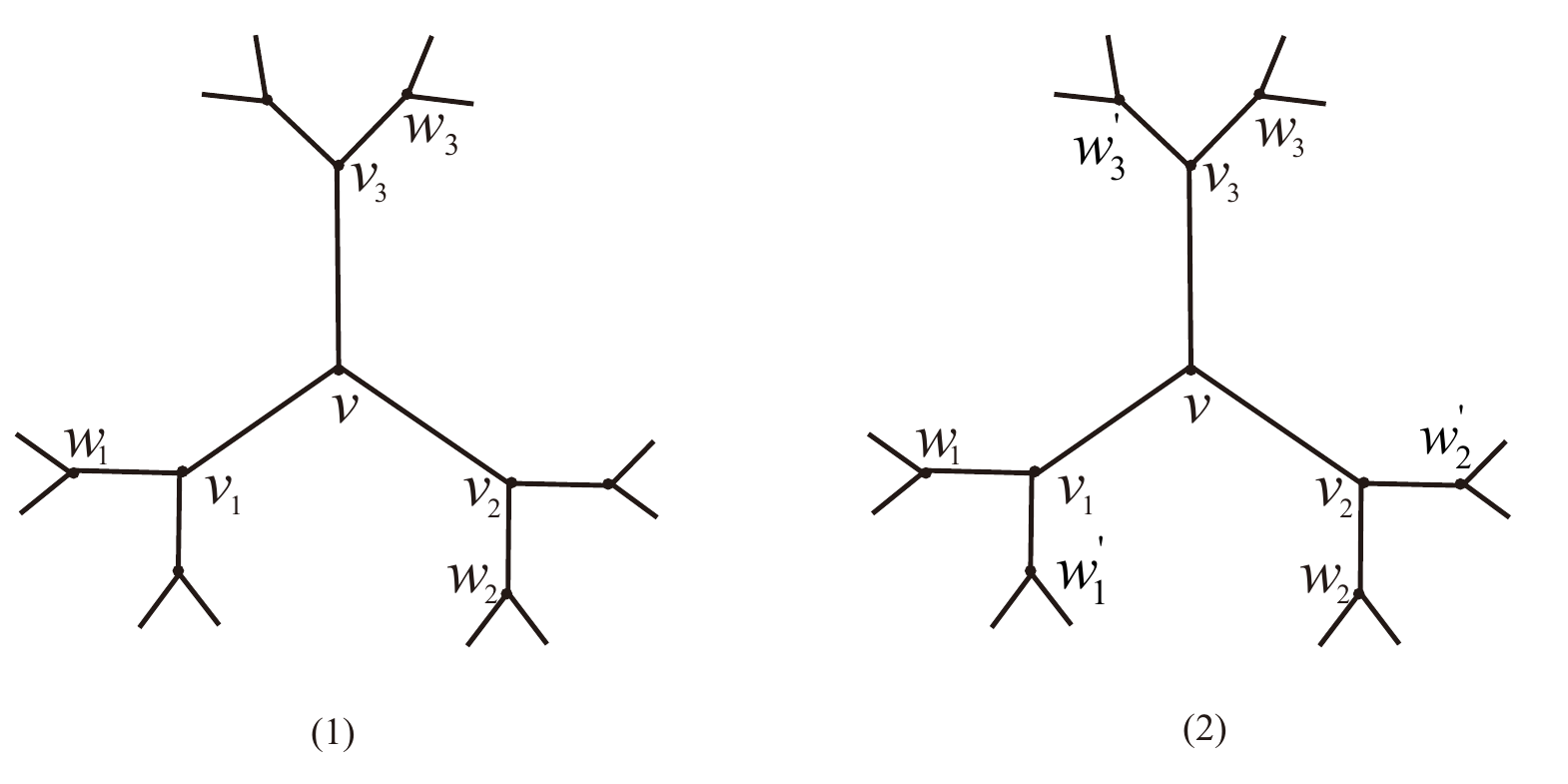}
\caption{}
\label{C}
\end{figure}

\begin{lem}\label{lemma4.1}
Each precoloring of $v_1w_1, v_2w_2, v_3w_3$ from their lists can be extended to an $L$-coloring of $H=G-v$.
\end{lem}

\pf 
 Order the edges in $G$ with respect to the distance from $v$, that is, if edge $e$ precedes edge $f$, then $d(v,e)\geq d(v,f)$, where $d(v,e)$ and $d(v,f)$ are the distance from $v$ to the edges $e$ and $f$, respectively. Then the last three edges in the list are $vv_1, vv_2$ and $vv_3$. The edges $v_iw_i$ and $v_iw_i'$ for $i\in [3]$ precede them. Color the edges in the list from the first to the last greedily. For each $e=xy$ in the list with $d(v,x)\ge d(v,y)\ge 1$, $y$ is adjacent to some vertex $z$ with $d(v,z)<d(v,y)$. So the three edges incident with $z$ are after the edge $e$ in the list. Clearly, at least two of the three edges at $z$ are not precolored, thus $e$ sees at least two uncolored edges in $G$ (Figure 3). So $e$ sees at most $10$ different colors, and thus can be colored. $\Box$
\begin{figure}[htbp]
\centering
\includegraphics[scale=0.5]{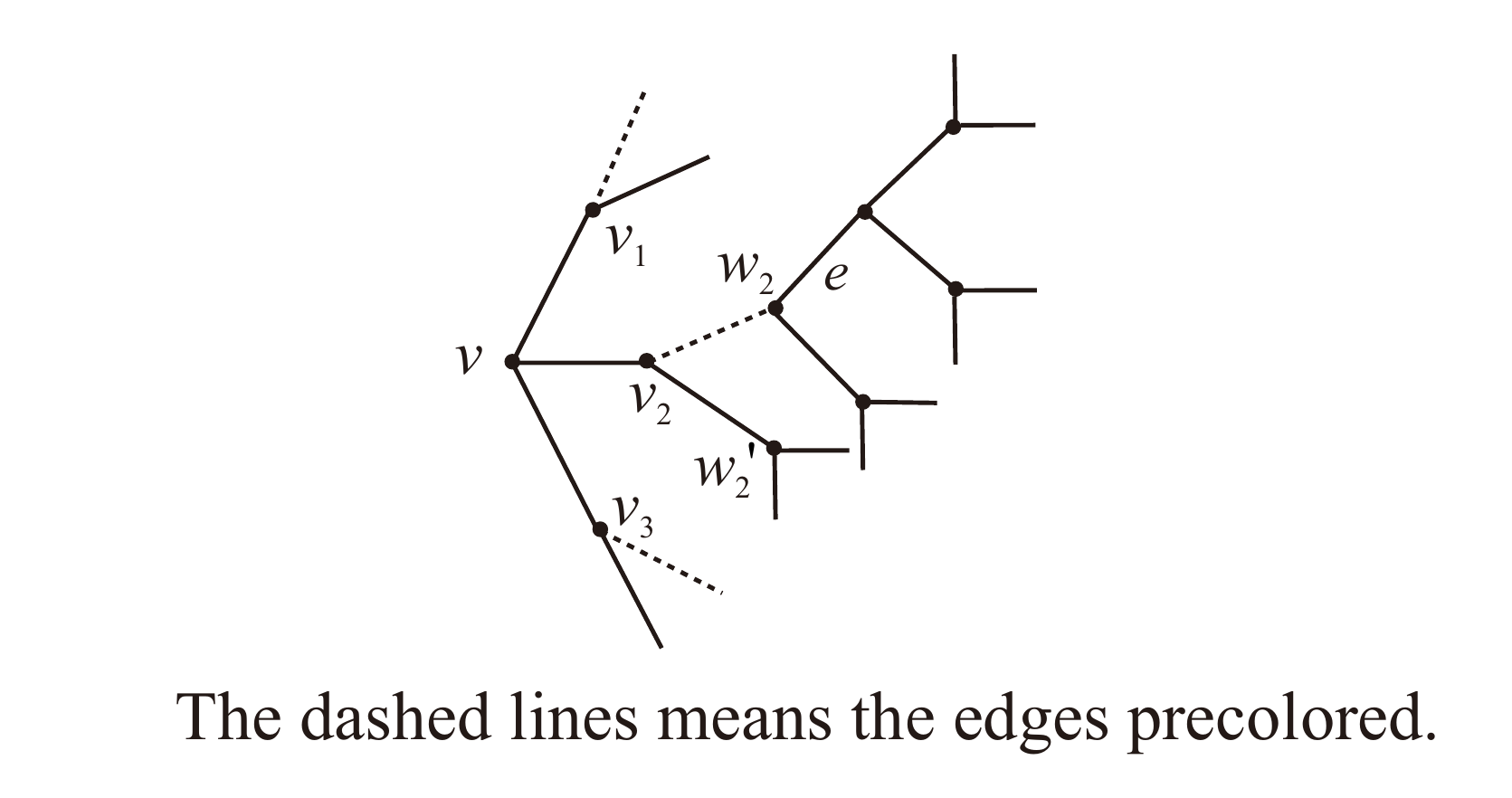}
\caption{}
\label{C}
\end{figure}
\medskip

For $i\in [3]$, let $B_i=L(v_iw_i)\cup L(v_iw_i')$ and $L_i=L(vv_i)$.  We will prove Theorem~\ref{th1} through a series of claims.

(1)  $B_1\cap B_2\cap B_3=\emptyset$.

For otherwise, precolor $v_iw_i$ (or $v_iw_i'$) with $\alpha\in B_1\cap B_2\cap B_3$, which can be extended to an $L$-coloring of $H$ by Lemma~\ref{lemma4.1}.  Now $vv_i$ for $i\in [3]$ sees at most $8$ different colors, so $|A(vv_i)|\ge 3$. So $vv_1, vv_2, vv_3$ can be colored in the order.

\medskip

(2)  For any $i,j\in [3]$ with $i\not=j$, $B_i\cap B_j\subseteq L_1\cap L_2\cap L_3$.

Suppose that for some $i, j\in [3]$ with $i\not=j$, there exists $\alpha\in (B_{i}\cap B_{j})-(L_1\cap L_2\cap L_3)$. It means $\alpha$ must not belong to the one of $L_1$, $L_2$, $L_3$. We assume $\alpha\in (B_{2}\cap B_{3})-L_1$. Without loss of generality, assume $\alpha\in L(v_{2}w_2)\cap L(v_{3}w_3)-L_1$. Precolor $v_2w_2, v_3w_3$ with $\alpha$ and by Lemma~\ref{lemma4.1}, we can extend it to an $L$-coloring $\phi$ of $H$. Then $|A_{\phi}(vv_1)|\geq3$, $|A_{\phi}(vv_2)|\geq2$, $|A_{\phi}(vv_3)|\geq2$, and we can color $vv_2, vv_3, vv_1$ by Lemma 1.4.

\medskip

(3) For some $i,j\in [3]$ with $i\not=j$, $L_i\cap L_j\not=\emptyset$.

For otherwise, in an $L$-coloring of $H$, each of $vv_1, vv_2, vv_3$ has an available color and the colors are distinct, so they could be colored.

\medskip
We may assume that $L_2\cap L_3\neq\emptyset$.\par
\medskip
(4)$|\bigcup_{i=1}^{3} B_{i}|\ge |\bigcup_{i=1}^{3}L_i|+|L_1\cap L_2\cap L_3|$.

By (1), $B_1\cap B_2, B_2\cap B_3, B_3\cap B_1$ are disjoint, and by (2), are subsets of $L_1\cap L_2\cap L_3$. So
\begin{align*}
|\bigcup_{i=1}^{3} B_{i}|&=\sum _{i=1}^{3}|B_{i}|-\sum_{i,j\in [3]}|B_i\cap B_j|+|\bigcap_{i=1}^{3}B_{i}|
\geq \sum _{i=1}^{3}|B_{i}|-|\bigcap_{i=1}^{3}L_i|=33-|\bigcap_{i=1}^{3}L_i|.
\end{align*}

On the other hand,  $$|\bigcup_{i=1}^{3}L_i|=\sum _{i=1}^{3}|L_i|-\sum_{i,j\in [3]}|L_i\bigcap L_j|+|\bigcap_{i=1}^{3}L_i|\leq 33-2|\bigcap_{i=1}^{3}L_i|.$$

Therefore, $|\bigcup_{i=1}^{3} B_{i}|\ge |\bigcup_{i=1}^{3}L_i|+|L_1\cap L_2\cap L_3|$.

\medskip

(5) For some $i,j\in [3]$ with $i\not=j$,  $B_{i}\bigcap B_{j}\not=\emptyset$.

For otherwise, $|B_{i}|\geq 11$ for $i=1,2,3$ and $|\bigcup_{i=1}^{3}B_i|\geq 33$. Since $L_2\cap L_3\neq\emptyset$, we have $|\bigcup_{i=1}^{3}L_i|\leq 32$. So there exists $\alpha\in (B_{1}\cup B_{2}\cup B_{3})-(L_1\cup L_2\cup L_3)$. Assume $\alpha\in L(v_1w_{1})\subset B_{1}$.  Since $|B_{2}\cup B_{3}|\geq 22$ and $|L_2\cup L_3|\leq 21$, there exists $\beta\in (B_{2}\cup B_{3})-(L_2\cup L_3)$, and we may assume $\beta\in L(v_2w_2)$.  Now we precolor $v_1w_1$ with $\alpha$ and $v_2w_2$ with $\beta$, and by Lemma~\ref{lemma4.1}, extend it to an $L$-coloring $\phi$ of $H$. Now $|A_{\phi}(vv_1)|\geq2$, $|A_{\phi}(vv_2)|\geq3$, $|A_{\phi}(vv_3)|\geq3$, we can color $vv_3, vv_2, vv_1$ by Lemma 1.4.

\bigskip

By (5) and (2), $|L_1\cap L_2\cap L_3|\ge 1$, so by (4), there  exists $\alpha\in \bigcup_{i=1}^{3} B_{i}-\bigcup_{i=1}^{3}L_i$. Assume $\alpha\in L(v_1w_1)\subseteq B_1$.

Precolor $v_1w_1$ with $\alpha$.

\begin{itemize}
\item $B_{2}\cap B_{3}\neq\emptyset$.

Let $\beta\in L(v_2w_2)\cap L(v_3w_3)$.  Precolor $v_2w_2, v_3w_3$ with $\beta$.  By Lemma~\ref{lemma4.1}, this precoloring can be extended to an $L$-coloring $\phi$ of $H$.  Note that for $i\in [3]$, $|A_{\phi}(vv_i)|\geq3$, we can color $vv_1, vv_2, vv_3$ in the order.

\item $B_{2}\cap B_{3}=\emptyset$.

Then $|B_{2}\cup B_{3}|=22-|B_{2}\cap B_{3}|=22>22-|L(e_{2})\cap L(e_{3})|=|L(e_{2})\cup L(e_{3})|$.  So there exists $\beta\in(B_{2}\cup B_{3})-(L_2\cup L_3)$. Suppose that $\beta\in L(v_3w_3)$ without loss of generality. Precolor $v_3w_3$ with $\beta$, by Lemma~\ref{lemma4.1}, this precoloring can be extended to an $L$-coloring $\phi$ of $H$.   Note that $|A_{\phi}(vv_1)|\geq2, |A_{\phi}(vv_2)|\geq3, |A_{\phi}(vv_3)|\geq3$, and we can color $vv_1, vv_2, vv_3$ in the order. $\Box$
\end{itemize}

\section{Final discussion}
As we mentioned in the introduction, one may try to improve our results by one, which, if true, would be optimal.  But this may not be easy, especially for subcubic planar graphs.

Here is another related question.  A graph is {\em chromatic-choosable} if its chromatic number equals to its list chromatic number.  It is an interesting problem to find graphs that are chromatic-choosable.  Zhu asked whether there exists a constant integer $k$ such that the $k$-th power $G^{k}$ is chromatic-choosable for every graph $G$. Kim, Kwon, and Park~\cite{KWP} answered this question negatively. Moreover, for any fixed $k$ they showed that there are graphs $G$ such that the value $\chi_{l}(G^{k})-\chi(G^{k})$ can be arbitrarily large.

We know $\chi'_{s,l}(G)$ is the list chromatic number of the square of the line graph of $G$. Kostochka and Woodall~\cite{AV} asked whether $G^2$ is chromatic-choosable for every graph. Kim and Park~\cite{KP} solved the conjecture in the negative by finding a family of graphs $G$ whose squares are complete multipartite graphs with partite sets of unbounded size.
\begin{ques}
Is $G^2$ chromatic-choosable for every line graph $G$?
\end{ques}

{\bf Acknowledgement:} We are very grateful for the careful reading and many helpful comments from the referees.



\section*{Appendix}

Note that if $P(x_1,x_2,\cdots,x_m)$ is a polynomial with $\deg (P)=n$, and $k_1,k_2,\cdots,k_m$ are non-negative integers with $\sum\limits_{i=1}^{m} k_i=n$. Let $c_P(x_1^{k_1}x_2^{k_2}\cdots x_m^{k_m})$ be the coefficient of monomial $x_1^{k_1}x_2^{k_2}\cdots x_m^{k_m}$ in $P$. Then
$$\frac{\partial^n P}{\partial x_1^{k_1}\partial x_2^{k_2}\cdots\partial x_m^{k_m}}=c_P(x_1^{k_1}x_2^{k_2}\cdots x_m^{k_m})\prod_{i=1}^{m}k_i!.$$

\lstset{
  language=matlab,
  basicstyle=\small,
  stepnumber=1, 
  numbersep=5pt,
commentstyle=\color{red!50!green!50!blue!50},
showstringspaces=false,
basicstyle=\footnotesize,
xleftmargin=2em,xrightmargin=2em,
  tabsize=4,
}
\begin{lstlisting}
%input
syms x1 x2 x3 x4 x5 x6 x7 x8 x9
%Lemma 2.4
Q=(x1-x2)*(x1-x3)*(x1-x4)*(x1-x5)*(x1-x6)*(x1-x7)*(x1-x8)*(x2-x3)*(x2-x4)*(x2-x5)*(x2-x6)*(x2-x7)*(x2-x8)*(x2-x9)*(x3-x4)*(x3-x5)*(x3-x7)*(x3-x8)*(x3-x9)*(x4-x5)*(x4-x6)*(x4-x8)*(x4-x9)*(x5-x6)*(x5-x7)*(x5-x9)*(x6-x7)*(x6-x8)*(x6-x9)*(x7-x8)*(x7-x9)*(x8-x9);
C=diff(diff(diff(diff(diff(diff(diff(diff(diff(Q,x1,4),x2,5),x3,4),x4,4),x5,3),x6,3),x7,3),x8,3),x9,3)/factorial(4)/factorial(5)/factorial(4)/factorial(4)/factorial(3)/factorial(3)/factorial(3)/factorial(3)/factorial(3)
%output
C=-6
\end{lstlisting}

\end{document}